\documentclass[a4paper,9pt]{article}

\usepackage{amsmath,amssymb}
\def\uqg{U_q(\hat{\mathfrak{gl}}_n)}

\def\uqgr{U_q^{\mathrm{res}}(\hat{\mathfrak{gl}}_n)}
\def\C{\mathbb{C}}
\def\Z{\mathbb{Z}}
\def\N{\mathbb{N}}
\def\V{\mathbf{V}}
\def\M{\mathbf{M}}
\def\R{\mathbf{R}}
\def\S{{\mathfrak{S}}}
\newtheorem{theo}{\bf{Theorem}}
\newtheorem{lem}{Lemma}
\newtheorem{cor}{Corollary}
\begin{document}
\title{Quantum affine algebras at roots of unity \\
and equivariant K-theory}
\author{Olivier Schiffmann}
\maketitle
\begin{abstract}
We show that the algebra homomorphism $U_q(\hat{\mathfrak{gl}}_n)
  \twoheadrightarrow K^{GL_d \times \C^*}(Z) \otimes \C(q)$
  constructed by Ginzburg and Vasserot between the quantum affine
  algebra of type $\mathfrak{gl}_n$ and the equivariant K-theory group
  of the Steinberg variety (of incomplete flags), specializes to a
  surjective homomorphism
$U_{\epsilon}^{\mathrm{res}}(\hat{\mathfrak{gl}}_n) \twoheadrightarrow
  K_\epsilon^{GL_d \times \C^*} (Z)$. In particular, this shows that
  the parametrization of irreducible
  $U_\epsilon^{\mathrm{res}}(\hat{\mathfrak{sl}}_n)$-modules and the
  multiplicity formulas of \cite{GV},\cite{V} are still valid when
  $\epsilon$ is a root of unity.
\end{abstract}
\vspace{.4in}
We will use the following notations: let $n,d \in \N$. We set $\V_d=\{(v_1,v_2
  \ldots v_n) \in \N^n\,|\,\sum_iv_i=d\}$; for $v \in \V_d$ put
  $\overline{v}_i=v_1 + \ldots + v_i$ and identify $v$ with the point
  $\sum_i v_i e_i \in \C^n$ where $(e_i)$ is the canonical basis of $\C^n$. For $v,w \in \V_d$, let
$$\mathbf{M}(v,w)=\big\{\,A=(a_{ij})_{i,j=1}^n\in \N^{n^2} \,|\, \sum_{j}a_{ij}=v_i,\;\sum_{i}a_{ij}=w_j\,\big\}$$
and $\M=\bigsqcup_{v,w}M(v,w)$. For $l \geq 1$, we denote by $\S^l$ 
the symmetric group in $l$ variables and we put $\S^{(v)}=\S^{v_1} \times \ldots
\times \S^{v_n} \subset \S^d$ for $v \in \V_d$ and
$\S^{(A)}=\S^{a_{11}} \times \S^{a_{21}} \times \ldots \times
\S^{a_{nn}} \subseteq \S^d$ for $A \in \M$. For all $v,w \in
\mathbf{V}_d$ the set $\mathbf{M}(v,w)$ can be identified with $\S^{(v)}
\backslash \S^d /\S^{(w)}$ (cf \cite{BLM}). We will write $\leq$ for
the Bruhat order on $\S^d$ and again $\leq$ for the induced order on
$\mathbf{M}(v,w)$. We set $\R=\C[x_1^{\pm 1}, \ldots x_d^{\pm
  1},q^{\pm 1}]$, $\R^{(v)}=\R^{\S^{(v)}}$ and
$\R^{(A)}=\R^{\S^{(A)}}$. Finally, $\overline{X}$ will denote the
Zariski closure of a subset $X$ of an algebraic variety.
\section{K-theory of the Steinberg variety}
Let $F$ be the variety of n-step flags in $\C^d$:
$$F=\big\{ (D_i)\,|\, 0=D_0 \subseteq D_1 \subseteq \ldots \subseteq D_n=\C^d\big\}$$
Then $F=\bigsqcup_{v \in \V_d} F_v$ where $F_v$ is the connected component
of $F$ consisiting of flags $(D_i)$ satisfying
$\mathrm{dim}\,(D_i/D_{i-1})=v_i$. The group $GL_d(\C)$ acts
diagonally on $F \times F=\bigsqcup_{v,w}F_v \times F_w$ and the
corresponding orbits are parametrized by $\mathbf{M}$: to $A \in \mathbf{M}(v,w)$ corresponds
$$O_A=\big\{\,((D_i),(D'_i))\,|\,\mathrm{dim}\,(D_i \cap D'_j)=\sum_{h \leq i, k \leq j} a_{hk}\,\big\} \subset F_v \times F_w$$
and we have $\overline{O_A}=\bigcup_{B \leq A} O_B$. For $A \in \M$,
denote by $Z_A=T^*_{O_A}(F \times F) \subset T^*F^2$ the conormal
bundle, set $Z_{v,w}=\bigcup_{A \in \mathbf{M}(v,w)} \overline{Z_A}$, and $Z=\bigcup_{v,w}
Z_{v,w}$ (the Steinberg variety). The group $G=GL_d \times \C^*$
acts on $Z$: $GL_d$ acts diagonally and $z \in \C^*$ acts by rescaling
by $z^{-2}$ along the fibers of $T^*F^2$. Let
$K^G(X)$ be the Grothendieck group of coherent $G$-equivariant sheaves
on a $G$-variety $X$ and let $[\mathcal{H}]$ denote the class of a
sheaf $\mathcal{H}$. Then $K^G(F_v) \simeq \R^{{(v)}}$,
$K^G(O_A)\simeq \R^{{(A)}}$. For all $A \in \M$, set $Z_{\leq
  A}=\bigcup_{B \leq A}Z_B$ and $Z_{<A}=\bigcup_{B < A}Z_B$. This
filtration induces a filtration on $K^G(Z)$, whose associated graded is
$\mathrm{gr}(K^G(Z))=\bigoplus_{A} K^G(Z_A) \simeq \bigoplus_A
\R^{(A)}$ (cf \cite{V}).\\
\hbox to1em{\hfill}The convolution product in equivariant K-theory
$\star: K^G(Z) \otimes K^G(Z) \to K^G(Z)$ equips $K^G(Z)$ with a
$\C[q,q^{-1}]$-associative algebra structure which is compatible with the filtration (cf \cite{CG}, chap.5).
\paragraph{}Fix $a \in \N$ and let $v\in \V_{d-a}$ be a composition of
$d-a$ in $n$ parts. Set $E_{ij}(v;a)=\mathrm{diag}(v)+ aE_{ij} \in
\M(v+ae_i, v+ae_j)$. The orbits $O_{E_{i,i\pm 1}(v;a)}$ are closed and
the projection $p_1: O_{E_{i,i\pm 1}(v;1)} \to F_{v+e_i}$ are smooth
and proper, with fibers isomorphic to ${\mathbb{P}}^{v_i}$. Denote by $\mathcal{O}_{p_1}(k)$ and
$T^*p_1$ the $k$th Serre twist and the cotangeant sheaf relative to
the fibers of $p_1$. Finally let $\mathcal{E}_{i,v,k}$
(resp. $\mathcal{F}_{i,v,k}$) be the pullback of
$\mathrm{Det}(T^*p_1) \otimes \mathcal{O}_{p_1}(k)$ on
$Z_{E_{i,i+1}(v;1)}$ (resp. on $Z_{E_{i+1,i}(v;1)}$) which we view as
$G$-equivariant sheaves on $Z$.
\section{The quantum affine algebra}
Let $\uqg$ be the quantum loop algebra of $\mathfrak{gl}_n$
over $\C(q)$, with generators
$E_{i,k}$, $F_{i,k}$, $K_j^{\pm 1}$, $K_{j,l}$ for $i=1,\ldots
n-1,\;j=1, \ldots n$, $k \in \Z$ and $j \in \Z^*$, with the relations of
Drinfeld's new presentation (cf \cite{DF}; we use the same notations
as in \cite{V}).\\
\hbox to1em{\hfill}For $\mathbf{\lambda}=\sum_{i=1}^n \lambda_i e_i
\in \Z^n$  let $\pi_{\mathbf{\lambda}}$ denote the projector on the
space of weight $\lambda$: the operators
$\pi_{\mathbf{\lambda}}$ satisfy the following relations: 
\begin{equation}\label{E:rel1}
\pi_{\mathbf{\lambda}}\pi_{\mathbf{\mu}}=\delta_{\mu,\lambda}\pi_{\mu}, \qquad K_i \pi_{\mathbf{\lambda}}=\pi_{\mathbf{\lambda}}K_i=K_{i,l}
\pi_{\mathbf{\lambda}}=\pi_{\mathbf{\lambda}}K_{i,l}=q^{\lambda_i}
\pi_{\mathbf{\lambda}},
\end{equation}
\begin{equation}\label{E:rel2}
E_{i,k}\pi_{\mathbf{\lambda}}=\pi_{\mathbf{\lambda}+e_i-e_{i+1}} E_{i,k}, \qquad
F_{i,k}\pi_{\mathbf{\lambda}}=\pi_{\mathbf{\lambda}-e_i+e_{i+1}}
F_{i,k}.
\end{equation}
By definition, the modified quantum loop algebra is
$\widetilde{U}_q(\hat{\mathfrak{gl}}_n)=\bigoplus_{\mathbf{\lambda}} \uqg
\pi_{\mathbf{\lambda}}$ (cf \cite{L}). The algebras $\uqg$ and
$\widetilde{U}_{q}(\hat{\mathfrak{gl}}_n)$ share the same
finite-dimensional representation theory.
\section{The map $\widetilde{U}_q(\hat{\mathfrak{gl}}_n) \twoheadrightarrow K^{GL_d \times \C^*}(Z) \otimes \C(q)$}
\begin{theo}[Ginzburg-Vasserot] The assignements
\begin{equation}\label{E:01}
E_{i,k}\pi_{v+e_{i+1}} \mapsto  (-q)^{-v_i} [\mathcal{E}_{i,v,k}], \qquad
F_{i,k}\pi_{v+e_{i}} \mapsto (-q)^{1-v_i} [\mathcal{F}_{i,v,k}] 
\end{equation}
for  $v \in \V_{d-1}$ and $E_{i,k}\pi_{v'}
\mapsto 0,\; F_{i,k}\pi_{v'} \mapsto 0$ otherwise, and
\begin{equation}\label{E:02}
\big(K_j^{\pm 1} + \sum_{l > 0} K_{i, \pm l}z^{\mp l}\big)\pi_v \mapsto
\prod_{m=1}^{\overline{v_i} -1} \frac{q^2z -x_m}{qz-qx_m}
\prod_{m=\overline{v_i}}^{d} \frac{z-q^2x_m}{qz-qx_m}
\in K^G(F_v)[[z^{\mp 1}]] \subset K^G(F)[[z^{\mp 1}]]
\end{equation}
for $v \in \mathbf{V}_d$, extend to a surjective algebra homomorphism
$\Phi:\widetilde{U}_q(\hat{\mathfrak{gl}}_n) \to K^G(Z)\otimes \C(q)$.
\end{theo}
In the last line, we identify $T^*F$ with the conormal bundle of the
diagonal of $F\times F$ in $T^*F^2$ and we consider $K^G(F)\simeq
K^G(T^*F)\subset K^G(Z)$. This theorem gives a geometric construction
of all irreducible, finite-dimensional
$U_q(\hat{\mathfrak{sl}}_n)$-modules when $q$ is not a root of unity.
\section{Restriction to integral forms}
\paragraph{}Let $\uqgr$ denote the restricted integral form of $\uqg$ (cf
\cite{CP2}): $\uqgr$ is the $\C[q,q^{-1}]$-algebra generated by
$E_{ik}^{(m)}=\frac{E_{ik}^m}{[m]!}$,
$F_{ik}^{(m)}=\frac{F_{ik}^m}{[m]!}$, $K_j^{\pm 1}$ and
$\frac{H_{j,r}}{[r]}$ where $H_{j,r}$ is determined by the relation
\begin{equation}\label{E:03}
\sum_{k \geq 1} H_{j, \pm k} z^{\mp k}= (q-q^{-1}) \mathrm{log}\, (1 + \sum_{l >0} K_{j}^{\mp 1}K_{j,\pm l} z^{\mp l})
\end{equation}
By definition, $\widetilde{U}_q^{\mathrm{res}}(\hat{\mathfrak{gl}}_n)=\bigoplus_\lambda \uqgr \pi_\lambda$.
\begin{theo} The map $\Phi: \widetilde{U}_q(\hat{\mathfrak{gl}}_n)
  \twoheadrightarrow K^G(Z) \otimes \C(q)$ restricts to a surjective map $\widetilde{\Phi}: \widetilde{U}_q^{\mathrm{res}}(\hat{\mathfrak{gl}}_n) \twoheadrightarrow K^G(Z)$.
\end{theo}
The theorem is a consequence of the following two lemmas.
\begin{lem}The map ${\Phi}$ restricts to a map $\tilde{\Phi}:\widetilde{U}_q^{\mathrm{res}}(\hat{\mathfrak{gl}}_n) \to K^G(Z)$. \end{lem}
\textbf{Proof:} it is clear that $\Phi(K_j^{\pm 1}\pi_v) \in K^G(Z)$
for all $v \in \mathbf{V}$; a direct computation shows that 
\begin{equation}\label{E:04}
\Phi\big(\frac{H_{j,r}}{[r]}\pi_v\big)=\frac{-1}{r} \Big(q^{\mp r} \sum_{l=1}^{\overline{v_i}-1} x_l^{\pm r} + q^{\pm r} \sum_{l=1+ \overline{v_i}}^d x_l^{\pm r}\Big) \in K^G(Z)
\end{equation}
For $v$ a composition of $d-1$ in $n$ parts, set $v^{(l)}=v +
l(e_{i+1}-e_i)$. Write $\overset{\to}{\underset{i=1,\ldots
    r}{\prod}}u_i$ for the ordered product $u_1u_2\ldots u_r$. We have
\begin{equation*}
\Phi (E_{i,k}^m\pi_{v^{(m-1)}+e_{i+1}})=\Phi(E_{i,k}\pi_{v^{(0)}+e_{i+1}})\star\ldots \star\Phi(E_{i,k}\pi_{v^{(m-1)}+e_{i+1}})=  \overset{\to}{\underset{l=0,\ldots m-1}{\prod}} (-q)^{v_i-l} [\mathcal{E}_{i,v^{(l)},k}].
\end{equation*}
 To see that $\Phi(E_{i,k}^{(m)}\pi_{v^{(m-1)}+e_{i+1}}) \in
 K^G(Z)$, it is enough to treat the case $n=2$. Let $v=(v_1,v_2)$,
 $v_1+v_2=d-1$. Recall that orbits of type $O_{E_{12}(w;a)}$
 are closed. Using \cite{V}, Lemme 12 et Exemple 4, we have 
\begin{equation}\label{E:45}
[\mathcal{E}_{1,v^{(l)},k}]= x_{1-l+v_1}^k \prod_{t=1}^{v_1-l}
\frac{x_t}{x_{1-l+v_1}} \in \R^{(E_{12}(v^{(l)};1))} =
K^G(Z_{E_{12}(v^{(l)};1)})
\end{equation}
and the convolution product
$$\star: K^G(Z_{E_{12}(v-le_1;l+1)}) \otimes K^G(Z_{E_{12}(v^{(l+1)};1)}) \to K^G(Z_{E_{12}(v-(l+1)e_1;l+2)})$$
 can be written
$$f \otimes g \to \S_{I \times J}^{I \cup J} \big(fg \prod _{i \in I}\prod_{j \in J} \frac{1-q^2x_j/x_i}{1-x_i/x_j}\big)$$
where $I=\{v_1-l\}$, $J=[v_1-l+1, v_1+1]$ and $\S_{I \times J}^{I \cup
  J}$ is the symmetrization map $\R^{\S_I \times \S_J} \to \R^{\S_{I
    \cup J}}$. It follows by induction that
\begin{equation}\label{E:05}
\overset{\to}{\underset{l=0 \ldots m-1}{\prod}} [\mathcal{E}_{1,v^{(l)},k}]=q^{\frac{m(m-1)}{2}}[m]! x_{2-m+v_1}^k \ldots x_{v_1+1}^k \prod_{t=1}^{1-m+v_1} \frac{x_t^m}{ x_{2-m+v_1} \ldots x_{1+v_1}}.
\end{equation}
Hence $\Phi(E_{i,k}^{(m)}\pi_{v^{(m-1)}+e_{i+1}}) \in K^G(Z)$. The
case of $F_{i,k}^{(m)}\pi_{\lambda}$ is similar.$\hfill\blacksquare$
\paragraph{}We now show that $\widetilde{\Phi}$ is
surjective. The algebra $K^G(Z)$ is generated by sheaves supported on
$\overline{Z_A}$, where $A \in \M$ is diagonal or of the type
$E_{i,i\pm 1}(v;a)$, $a \in \N$ (cf. \cite{V}, proof of
proposition 11). From (\ref{E:04}), we deduce that
$K^G(\overline{Z_A}) \subset \mathrm{Im}\;\widetilde{\Phi}$ for all
diagonal $A$. It thus remains to show the following:
\begin{lem} We have $K^G(\overline{Z_A}) \subset \mathrm{Im}\;\widetilde{\Phi}$ for $A=E_{i, i \pm 1}(v;a)$. \end{lem}
\textbf{Proof:} We treat the case $A=E_{i,i+1}(v;a)$ and proceed by
induction on $a$. We have
$K^G(Z_{E_{i,i+1}(v;a)})\simeq \R^{\S_{v_1} \times \ldots \times
  \S_{v_i} \times \S_a \times \S_{v_{i+1}} \ldots \times
  \S_{v_n}}$. The algebra $K^G(Z_{\mathrm{diag}(v+ae_i)})=\R^{\S_{v_1}
  \times \ldots \times \S_{a+v_i} \ldots \times \S_{v_n}}$ acts (by
convolution) on $K^G(Z_{E_{i,i+1}(v;a)})$. For $a=1$, the elements
$[\mathcal{E}_{i,v,k}]$, $k \in \Z$ form a generating system of the $K^G(Z_{\mathrm{diag}(v+e_i)})$-module $K^G(Z_{E_{i,i+1}(v;1)})$
(cf (\ref{E:45})). Now let $a>1$ and suppose that
$K^G(Z_{E_{i,i+1}(v';b)}) \subset \mathrm{Im}\,\widetilde{\Phi}$ for
all $v'$ and all $b <a$. The $K^G(Z_{\mathrm{diag}(v+ae_i)})$-module
$K^G(Z_{E_{i,i+1}(v;a)})$ is generated by the subspace $\C[y_1^{\pm
  1},
\ldots y_a^{\pm 1}]^{\S^a}$ where we set
$y_l=x_{l+\overline{v}_i}$. Let us denote by
$\S:\C[y_1^{\pm 1},\ldots y_a^{\pm 1}] \to \C[y_1^{\pm 1},\ldots
y_a^{\pm 1}]^{\S_a}$ the operator of complete symmetrization, so that
we have $\C[y_1^{\pm 1},
\ldots y_a^{\pm 1}]^{\S^a}=\bigoplus_\alpha
\C\S(y_1^{\alpha_1}\ldots y_a^{\alpha_a})$ for $\alpha=(\alpha_i) \in
\Z^a$. We show, by induction on $\parallel \alpha \parallel =\sum_t
(\alpha_t-\mathrm{inf}\,\alpha)^2$, that $\S(y_1^{\alpha_1}\ldots
y_a^{\alpha_a}) \in \mathrm{Im}\;\widetilde{\Phi}$. Using (\ref{E:05}) with
$m=a$ we obtain
$$\widetilde{\Phi}(E_{i,l}^{(a)}\pi_{v^{(a-1)}+(a-1)e_{i}+e_{i+1}})=(-1)^aq^{a(a-1+v_i)}y_1^l\ldots
y_a^l \prod_{t=1}^{\overline{v}_i} \frac{x_t^a}{y_1\ldots y_a}.$$
It follows that $y_1^k\ldots y_a^k \in \sum_l
K^G(Z_{\mathrm{diag}(v+ae_{i+1})})
\widetilde{\Phi}(E_{i,l}^{(a)}\pi_{v^{(a-1)}+(a-1)e_{i}+e_{i+1}})
\subset \mathrm{Im}\,\widetilde{\Phi}$ for all $k \in \Z$. Now let
$h>0$ and suppose that $\S(y_1^{\alpha_1}\ldots y_a^{\alpha_a}) \in
\mathrm{Im}\,\widetilde{\Phi}$ for all $\alpha$ such that
$\parallel \alpha \parallel <h$. Fix some $\alpha=(\alpha_1=\ldots
=\alpha_s>\alpha_{s+1}\geq\ldots \geq \alpha_a)$, $\parallel
\alpha \parallel =h$. Denote by $\S':\C[y_{s+1}^{\pm 1},\ldots y_{a}^{\pm
  1}]\to
\C[y_{s+1}^{\pm 1},\ldots y_{a}^{\pm 1}]^{\S_{a-s}}$ the operator of
symmetrization and set
 \begin{align*}
P&=\S'(y_{s+1}^{\alpha_{s+1}}\ldots y_a^{\alpha_a})\prod_{u=s+1}^a y_u
\in \C[y_{s+1}^{\pm 1},
\ldots y_a^{\pm 1}]^{\S^{a-s}} \subset K^G(Z_{E_{i,i+1}(v+se_i;a-s)}),\\
Q&=y_1^{\alpha_1} \ldots y_s^{\alpha_s}\prod_{t=1}^s\frac{1}{y_t}
\in \C[y_1^{\pm 1},\ldots y_s^{\pm 1}]^{\S^s} \subset
K^G(Z_{E_{i,i+1}(v+(a-s)e_{i+1};s)}).
\end{align*}
 Then
\begin{equation}\label{E:6}
P \star Q = \S\Big(y_1^{\alpha_1} \ldots
y_s^{\alpha_s}\S'(y_{s+1}^{\alpha_{s+1}} \ldots y_a^{\alpha_a})
\prod_{t=1}^ s\prod_{u=s+1}^a (1+(1-q^2)\frac{y_u}{y_t-y_u})\Big)
\end{equation}
Using the relation
$(y_t^{\alpha}y_u^{\beta+1}-y_t^{\beta+1}y_u^{\alpha})/(y_t-y_u)=y_t^{\alpha-1}y_u^{\beta+1}+
\ldots+ y_t^{\beta+1}y_u^{\alpha-1}$, and noticing that $(\alpha-i)^2 +
(\beta +i)^2 < \alpha^2 + \beta^2$ if $\alpha > \beta\geq 0$ and $i<
\alpha-\beta$, it is easy to see that the r.h.s of (\ref{E:6}) is
equal to $ (a-s)!\S(y_1^{\alpha_1} \ldots y_a^{\alpha_a}) + T$
where $T$ is a linear combination of polynomials
$\S(y_1^{\beta_1}\ldots y_a^{\beta_a})$ with $\sigma_2(\beta) <
\sigma_2(\alpha)=h$. Since $s,a-s <a$, it follows by the induction
hypothesis that
$\S(y_1^{\alpha_1}\ldots y_a^{\alpha_a}) \in
\mathrm{Im}\,\widetilde{\Phi}$. Hence, $\C[y_1^{\pm 1},\ldots y_a^{\pm
  1}]^{\S^a}
\subset \mathrm{Im}\,\widetilde{\Phi}$ and $K^G(Z_{E_{i,i+1}(v;a)})
\subset \mathrm{Im}\,\widetilde{\Phi}$. This concludes the induction
and the proof of the theorem. $\hfill\blacksquare$

\paragraph{}There is a surjective morphism
$\Psi:U_q(\hat{\mathfrak{gl}}_n) \twoheadrightarrow K^G(Z) \otimes
  \C(q)$ analogous to $\Phi$, (defined by $\Psi(u)=\sum_{\lambda \in
  \V_d}\Phi(u\pi_\lambda)$ for $u \in
  U_q(\hat{\mathfrak{gl}}_n)$). For $\epsilon \in \C^*$, let us set
  $U_\epsilon^{\mathrm{res}}(\hat{\mathfrak{gl}_n})=U_q^{\mathrm{res}}(\hat{\mathfrak{gl}}_n)_{|q=\epsilon}$ and $K_\epsilon^G(Z)=K^G(Z)_{|q=\epsilon}$. Lemma 1 shows that $\Psi$ restricts to a map $\widetilde{\Psi}:U_q^{\mathrm{res}}(\hat{\mathfrak{gl}}_n) \to K^G(Z)$, and specializes to a morphism $\Psi_\epsilon: U_\epsilon^{\mathrm{res}}(\hat{\mathfrak{gl}}_n) \to K^G_\epsilon(Z)$.
\begin{cor} For all $\epsilon \in \C^*$, the map $\Psi_\epsilon: U_\epsilon^{\mathrm{res}}(\hat{\mathfrak{gl}}_n) \to K^G_\epsilon(Z)$ is surjective.
\end{cor}
Indeed, for all $\epsilon \in \C^*$, the projectors on weight
subspaces $\pi_\lambda$ can be realized inside $U_\epsilon^{\mathrm{res}}(\hat{\mathfrak{gl}}_n)$, i.e
there exists $\overline{\pi}_\lambda\in
U_\epsilon^{\mathrm{res}}(\hat{\mathfrak{gl}}_n)$ satisfying
(\ref{E:rel1}-\ref{E:rel2}). Let us denote by
$i_\epsilon:\widetilde{U}_\epsilon^{\mathrm{res}}(\hat{\mathfrak{gl}}_n)\to
U_\epsilon^{\mathrm{res}}(\hat{\mathfrak{gl}}_n)$
the algebra map defined by $i_\epsilon(u \pi_\lambda)=u
\overline{\pi}_\lambda$ for $u \in
U_\epsilon^{\mathrm{res}}(\hat{\mathfrak{gl}}_n)$. We have
$\widetilde{\Phi_\epsilon}=\widetilde{\Psi} \circ i_\epsilon$ and the
surjectivity of $\widetilde{\Psi}_\epsilon$ now follows from the
surjectivity of
$\widetilde{\Phi}_\epsilon$.$\hfill\blacksquare$
\paragraph{}Theorem 2 justifies the constructions of
$U_\epsilon^{\mathrm{res}}(\hat{\mathfrak{sl}}_n)$-modules and the
Kazhdan-Lusztig multiplicity formulas of \cite{GV},\cite{V} when
$\epsilon$ is a root of unity:
\begin{cor} The parametrization of irreducible, finite-dimensional
  $U_\epsilon^{\mathrm{res}}(\hat{\mathfrak{sl}}_n)$-modules $L_O$
  (resp. standard modules $M_O$) in terms of graded nilpotent orbits
  $O$ (cf \cite{V}) and the multiplicity formula 
$$[M_O:L_{O'}]=\sum_i \mathrm{dim}\,\mathcal{H}^i(IC_{O'})_{|O}$$
are valid when $\epsilon$ is a root of unity.\end{cor}

\paragraph{Acknowledgments:} I would like to thank Eric Vasserot for
his patience and precious advice.

\small{
\vspace{4mm}
Olivier Schiffmann, ENS Paris, 45 rue d'Ulm, 75005
PARIS; \texttt{schiffma@clipper.ens.fr}
\end{document}